\documentclass{amsart}
\usepackage{amsxtra}
\theoremstyle{plain}

\newcommand{\cleqn}{\setcounter{equation}{0}}
\newcommand{\clth}{\setcounter{theorem}{0}}
\newcommand {\sectionnew}[1]{\section{#1}\cleqn\clth}
\newcommand{\nn}{\hfill\nonumber}
\newtheorem{theorem}{Theorem}[section]
\newtheorem{lemma}[theorem]{Lemma}
\newtheorem{definition-theorem}[theorem]{Definition-Theorem}
\newtheorem{proposition}[theorem]{Proposition}
\newtheorem{corollary}[theorem]{Corollary}
\newtheorem{definition}[theorem]{Definition}
\newtheorem{example}[theorem]{Example}
\newtheorem{remark}[theorem]{Remark}
\newcommand \bth[1] { \begin{theorem}\label{t#1} }
\newcommand \ble[1] { \begin{lemma}\label{l#1} }

\newcommand \bpr[1] { \begin{proposition}\label{p#1} }
\newcommand \bco[1] { \begin{corollary}\label{c#1} }
\newcommand \bde[1] { \begin{definition}\label{d#1}\rm }
\newcommand \bex[1] { \begin{example}\label{e#1}\rm }
\newcommand \bre[1] { \begin{remark}\label{r#1}\rm }

\newcommand {\eth} { \end{theorem} }
\newcommand {\ele} { \end{lemma} }

\newcommand {\epr} { \end{proposition} }
\newcommand {\eco} { \end{corollary} }
\newcommand {\ede} { \end{definition} }
\newcommand {\eex} { \end{example} }
\newcommand {\ere} { \end{remark} }

\newcommand \thref[1]{Theorem \ref{t#1}}
\newcommand \leref[1]{Lemma \ref{l#1}}
\newcommand \prref[1]{Proposition \ref{p#1}}
\newcommand \coref[1]{Corollary \ref{c#1}}

\newcommand \exref[1]{Example \ref{e#1}}

\newcommand \lb[1]{\label{#1}}
\def \Rset {{\mathbb R}}         
\def \Cset {{\mathbb C}}
\def \Zset {{\mathbb Z}}

\def \Qset {{\mathbb Q}}
\def \Pset {{\mathbb P}}
\def \Rg   {\Rset_{\geq 0}}
\def \cp   {{\Cset \Pset^1}}
\def \g  {\mathfrak{g}}   
\def \gt  {\tilde{\mathfrak{g}}}   
\def \k  {\mathfrak{k}}
\def \b  {\mathfrak{b}}  
\def \f  {\mathfrak{f}}  
  
\def \h  {\mathfrak{h}}

\def \C  {{\mathcal{C}}}           
\newcommand \Aff  {\mathcal{AFF}}
\def \Fin  {\mathcal{F} in}

\def \HC  {{\mathcal{HC}}}
\def \GF {{\mathcal{C}_{(\g, \f)}}}

\def \O  {{\mathcal{O}}}
\def \N  {{\mathcal{N}}}
\def \M  {{\mathcal{M}}}
\def \c  {{K}}

\def \de {\delta}
\def \al {\alpha}

\def \Ga {\Gamma}
\def \ga {\gamma}
\def \ka {\kappa}
\def \Om {\Omega}
\def \om {\omega}
\def \e  {t}   
\def \la {\lambda}

\def \vp {\varphi}
\def \sig{\sigma}

\def \op {\oplus}
\def \mt  {\mapsto}
\def \lla {\longleftarrow}
\def \ra  {\rightarrow}           

\def \hra {\hookrightarrow}
\def \sub {\subset}
\def \sup {\supset}
\def \st  {\ast}                 

\def \ol {\overline}            

\def \wh {\widehat}

\def \id { {\mathrm{id}} }
\def \Ann { {\mathrm{Ann}} }

\DeclareMathOperator \Span { {\mathrm{span}} }  
\DeclareMathOperator \ad { {\mathrm{ad}} }

\DeclareMathOperator \Hom { {\mathrm{Hom}} }

\DeclareMathOperator \I { {\mathrm{I}} }
\DeclareMathOperator \Ind { {\mathrm{Ind}} }
\DeclareMathOperator \Irr { {\mathrm{Irr}} }
\DeclareMathOperator \Res { {\mathrm{Res}} }
\renewcommand \Im { {\mathrm{Im}} \,}
\renewcommand \Re { {\mathrm{Re}} \,}
\renewcommand \max { {\mathrm{max}} }
\newcommand \ev { {\mathrm{Ex}} }
\begin{document}
\title[Finiteness of the Kazhdan--Lusztig tensor product]
{
Categories of modules over an affine Kac--Moody algebra 
and finiteness of the Kazhdan--Lusztig tensor product
}
\author[Milen Yakimov]{Milen Yakimov}
\address{
Department of Mathematics \\
University of California \\
Santa Barbara, CA 93106, U.S.A.}
\email{yakimov@math.ucsb.edu}
\date{}
\begin{abstract}
To each category $\C$ of modules of finite length over a complex simple 
Lie algebra $\g,$ closed under tensoring with finite dimensional 
modules, we associate and study a category $\Aff(\C)_\ka$ of smooth 
modules (in the sense of Kazhdan and Lusztig \cite{KL}) 
of finite length over the corresponding affine Kac--Moody algebra in
the case of central charge less than the critical level. 
Equivalent characterizations of these categories are obtained in the 
spirit of the works of Kazhdan--Lusztig \cite{KL}
and Lian--Zuckerman \cite{LZ, LZ2}. In the main part of 
this paper we establish a finiteness result for the 
Kazhdan--Lusztig tensor product which can be considered 
as an affine version of a theorem of Kostant \cite{Kos}.
It contains as special cases the finiteness results of  
Kazhdan, Lusztig \cite{KL} and Finkelberg \cite{F}, 
and states
that for any subalgebra $\f$ of $\g$ which is reductive in $\g$ 
the ``affinization'' of the category of finite length
admissible $(\g, \f)$ modules is stable under 
Kazhdan--Lusztig's tensoring
with the ``affinization'' of the category of finite dimensional 
$\g$ modules (which is $\O_\ka$ in the notation of 
\cite{KL, KL2, KL4}). 
\end{abstract}
\maketitle
\sectionnew{Introduction}\lb{Intro}
Let $\g$ be a complex simple Lie algebra
and let $\gt$ be the corresponding untwisted 
affine Kac--Moody algebra, 
see \cite{Kac}. It is the central extension of the loop algebra 
$\g[\e, \e^{-1}] = \g \otimes_\Cset \Cset[\e, \e^{-1}]$ by
\begin{equation}
[x \e^n, y \e^m]= [x,y] \e^{n+m} + n \delta_{n, -m}(x, y) \c, \; x, 
y \in \g
\lb{cocyc}
\end{equation}
where $(.,.)$ denotes the invariant bilinear form on $\g$ 
normalized by $(\al, \al)=2$ for long roots $\al.$ 
The affine Kac--Moody algebra $\gt$ 
is a $\Zset$ graded Lie algebra by 
\begin{equation}
\lb{grading}
\deg x \e^n = n, \; \deg \c =0.
\end{equation}
Set for shortness $\gt_+=\e \g[\e] \hra \gt.$ 
Denote the graded components of $U(\gt)$ and $U(\gt_+)$ of degree $N$ 
by $U(\gt)^N$ and $U(\gt_+)^N,$ respectively.

\bde{smooth}(Kazhdan--Lusztig) For a $\gt$ module $V$ define 
\[
V(N) = \Ann_{U(\gt_+)^N} V \sub V, \; N \in \Zset_{>0}.
\]
A $\gt$ module $V$ is called {\em{strictly smooth}} if
\[
\cup_N V(N) = V.
\]
\ede
Recall also that a $\gt$ module $V$ is called {\em{smooth}} 
if each vector in $V$ is annihilated 
by $\e^N \g[t]$ for all sufficiently large integers $N.$

Clearly 
\[
V(1) \sub V(2) \sub \dots
\]
and each $V(N)$ is a $\g$ module since $\g \hra \gt$ 
normalizes $U(\gt_+)^N.$ Each strictly smooth $\gt$ module is a module 
for the topological algebra $\hat{\g}$ which is the 
central extension of $\g((\e))$ by \eqref{cocyc}. 

The main objects of our consideration are the following categories
of $\gt$ modules.

\bde{main}
Let $\C$ be a full subcategory of the category of $\g$ modules
and $\ka \in \Cset.$
Define $\Aff(\C)_\ka$ to be the full subcategory of the category 
$\gt$ modules of central charge $\ka-h\spcheck$ consisting of 
strictly smooth,
finitely generated $\gt$ modules $V$ such that  
\[
V(N) \in \C, \; \mbox{for all} \; N=1, 2, \dots
\]
\ede
\noindent
Here $h\spcheck$ denotes the dual Coxeter number of $\g.$ 

We will restrict our attention only to categories $\C$  
of finite length $\g$ modules which are closed 
under tensoring with the adjoint representation and 
taking subquotients. The following are some important 
examples.  

\bex{mainex} 

(1) Let $\C =\Fin_\g$ be the category of finite dimensional 
$\g$ modules. Then $\Aff(\Fin_g)_\ka$ is Kazhdan--Lusztig's category
$\O_\ka$ defined in \cite{KL}.

(2) Let $\f$ be a subalgebra of $\g$ which is reductive in $\g,$ i.e.
$\g$ is completely reducible as an $\f$ module under the adjoint
representation, see \cite[section 1.7]{Dix} for details. 
Consider the 
category of finite length admissible $(\g, \f)$ modules -- 
$\g$ modules which restricted to $\f$ 
decompose to a direct sum of finite dimensional irreducible $\f$ 
modules, each of which occurs 
with finite multiplicity. It will be denoted
by $\GF.$ Kostant's theorem states that $\GF$ is  
closed under tensoring with finite dimensional $\g$ 
modules, \cite[Theorem 3.5]{Kos}. 
Clearly $\GF$ is also closed under taking subquotients. 
The example (1) is obtained from (2) 
when one specializes $\f=\g.$

(3) Let $\g_0$ be a real form of $\g$ and $\k$ be the complexification of 
a maximal compact subalgebra of $\g_0.$ Then $\k$ is reductive in $\g$ and 
the category from part 2 specializes to the category of Harish-Chandra 
$(\g, \k)$ modules to be denoted by $\HC_{(\g, \k)}.$ 
 
Lian and Zuckerman \cite{LZ, LZ2} studied the category of
$\Zset$-graded $\gt$ modules 
\begin{equation}
\label{LZ}
V=\oplus_{n \in \Zset} V_n
\end{equation}
for the grading \eqref{grading} of $\gt$
for which $V_n$ is a Harish-Chandra $\g$-module and $V_n=0$ for
$n \gg 0.$ Later we will show that this is exactly the category 
$\Aff(\HC_{(\g, \k)})$ (see \prref{LZchar} below) 
and will obtain a characterization
of all categories $\Aff(\C)_\ka$ in this spirit. 
At the same time the above $\Zset$-grading
is not canonical and is not preserved in general by $\gt$ 
homomorphisms. In the general case there is a canonical
$\Cset$ grading with similar properties, obtained from
the generalized eigenvalues of the Sugawara operator $L_0,$
see \prref{Sugchar}. It is preserved by arbitrary $\gt$ 
homomorphisms.

I. Frenkel and Malikov \cite{FM} defined a category 
of affine Harish-Chandra bimodules from the point of view of the  
construction of Bernstein and S. Gelfand \cite{BG}, treating 
case when $\g_0$ is a complex simple Lie algebra considered 
as a real algebra. It is unclear how in this case
\cite{FM} is related to \cite{LZ} and the constructions
of this paper. 

(4) Finally as another specialization of (2) one can choose $\f$ to be 
a Cartan subalgebra $\h$ of $\g$ and obtain the category of weight modules 
for $\g.$ Its subcategory $\O$ of Bernstein--Gelfand--Gelfand 
(for a fixed Borel subalgebra $\b \sup \h$ of $\g$) is also closed 
under tensoring with finite dimensional $\g$ modules and 
taking subquotients. From the results in Section 3, in particular
\thref{affineC}, it follows that the category $\Aff(\O)_\ka$ for
$\ka \notin \Rset_{\geq 0}$ 
is essentially the affine BGG category $\O$ of $\gt$ modules 
with central charge $\ka - h\spcheck,$ see \cite{Kac}. 
The generator $d$ for 
the extended affine Kac--Moody algebra acts on $\gt$ modules from 
$\Aff(\O)_\ka$ by $const-L_0$ where $L_0$ is the 0th Sugawara operator, 
see \cite{Kac} and Section 2 below. 
\eex

The following Theorem summarizes some of the main properties 
of the categories $\Aff(\C)_\ka.$    

\bth{mainth} Assume that the category $\C$ of finite length $\g$ 
modules is closed under tensoring with the adjoint representation and 
taking subquotients, and that 
$\ka \notin \Rset_{\geq 0}.$ Then the following hold:
 
(1) For any $M \in \C$ the induced module 
\[
\Ind(M)_\ka = U(\gt) \otimes_{U(\g[\e] \oplus \Cset \c)} M
\]
belongs to $\Aff(\C)_\ka.$ (In the definition of the tensor product 
the central element $\c$ acts on $M$ by
$(\ka - h\spcheck)$ and $\gt_+=\e \g[\e]$ annihilates $M.)$  

(2) The $\gt$ modules in $\Aff(\C)_\ka$ have finite length and 
are exactly the quotients of the
induced modules from $U(\g[\e])$ modules $\N$ annihilated by 
the degree $n$ component $U(\gt_+)^n$ of $U(\gt_+)$ for some 
$n \gg 0,$ recall \eqref{grading}.

(3) $\Aff(\C)_\ka$ is closed under taking subquotients.

(4) If the original category of $\g$ modules $\C$ 
is closed under extension then the category $\Aff(\C)_\ka$ is closed under 
extension inside the category of $\gt$ modules of 
central charge $\ka -h\spcheck.$ 
\eth
Let us also note that every irreducible module in 
$\Aff(\C)_\ka$ is (the) unique irreducible quotient of
$\Ind(M)_\ka$ for some irreducible $\g$ module $M.$ In addition 
for two nonisomorphic irreducible $\g$ modules the related 
irreducible $\gt$ modules are nonisomorphic.

Each smooth $\gt$ module of fixed central charge,
different from the critical level $-h\spcheck,$
canonically gives rise to a representation of the Virasoro
algebra by the Sugawara operators $L_k,$
see \cite{Kac} or the review in Section 2.
The generalized eigenspaces $V^\xi$ of the 
operator $L_0$ $(\xi \in \Cset)$ are naturally
$\g$ modules. One has the following characterization 
of $\Aff(\C)_\ka.$

\bpr{Sugchar} In the setting of \thref{mainth} the category 
$\Aff(\C)_\ka$ consists exactly of those 
finitely generated smooth $\gt$ modules 
of central charge $\ka - h\spcheck$
for which
\begin{equation}
\label{Sugdecom}
V = \bigoplus_{\xi \colon \xi -\xi_i \in \Zset_{\geq 0}} V^\xi \quad
\mbox{for some} \quad \xi_1, \; \ldots, \; \xi_n \in \Cset
\end{equation}
and $V^\xi \in \C.$
\epr

Kazhdan and Lusztig \cite{KL, KL2, KL4} defined a fusion tensor product
$V_1 \dot{\otimes} V_2$
of any two strictly smooth $\gt$ modules $V_1$ and $V_2,$
motivated by developments in conformal field theory
\cite{BPZ, KZ, MS, TUY, BFM}. In a related series of works
Huang and Lepowsky developed a theory of tensor products 
for modules over vertex operator algebras.
The modules $V_1 \dot{\otimes} V_2$ obtained from the Kazhdan--Lusztig 
tensor product are strictly smooth but in general it is hard to check 
under what conditions they have finite length. Kazhdan and Lusztig
proved in \cite{KL} that the category $\O_\ka$ is closed 
under the fusion tensor product and further used it to construct
functors to representations of quantized universal enveloping algebras.
In the case of positive integral central charge Cherednik
defined in \cite{C} a version of the fusion product for a category of 
integrable modules and showed that the latter is invariant under that 
tensor product.

For the affine algebra $\gt$ and the fusion tensor product 
the category $\O_\ka$ plays the role of the category of 
finite dimensional modules for the algebra $\g.$ {\em{We prove the 
following finiteness property of the Kazhdan--Lusztig tensor
product. It is an affine version of Kostant's theorem \cite{KL}
that tensoring with finite dimensional $\g$ modules preserves
the category of $\g$ modules $\GF.$}}

\bth{fin} Let $\ka \notin \Rset_{\geq 0}$ and
$\f$ be a subalgebra of the complex 
simple Lie algebra $\g$ which is reductive in $\g.$ 
Then the Kazhdan--Lusztig fusion tensor product 
of a $\gt$ module in $\O_\ka$ and a $\gt$ module in 
$\Aff(\GF)_\ka$ belongs to $\Aff(\GF)_\ka:$
\[
\dot{\otimes} \colon
\O_\ka \times \Aff(\GF)_\ka \ra \Aff(\GF)_\ka,
\quad
\dot{\otimes} \colon
\Aff(\GF)_\ka \times \O_\ka \ra \Aff(\GF)_\ka. 
\]
\eth

Thus $\Aff(\GF)_\ka$ become (bi)module categories for the ring 
category $\O_\ka$ using the left and right Kazhdan--Lusztig 
tensoring with objects from $\O_\ka.$ The associativity
is defined by a straightforward generalization of \cite{KL2}.
There are also braiding isomorphisms, intertwining the two tensor
products, again defined as in \cite{KL2}. 

In the special case $\f=\g$ one has 
$\O_\ka = \Aff(\C_{(\g,\g)})_\ka$
and we get just another proof of 
one of the main results of Kazhdan and Lusztig in \cite{KL} that
\[
\dot{\otimes} \colon
\O_\ka \times \O_\ka \ra \O_\ka.
\]
The novelty in this paper is a direct proof of \thref{fin}
which in the main part is independent of the 
one of Kazhdan and Lusztig \cite[Section 3]{KL} who use 
Soergel's generalized Bernstein--Gelfand--Gelfand (Brauer) reciprocity
\cite[Section 3.2]{BGS}.
They show that if $V$ is a strictly smooth $\gt$ module of central charge 
less than the critical level and $V(1)$ is finite dimensional then 
$V \in \O_\ka.$ Unfortunately BGG reciprocity does not 
to hold in the categories of $\gt$ modules $\Aff(\GF)_\ka$ and 
even in the standard category of Harish-Chandra $\g$ modules. (This was 
communicated to us by G. Zuckerman.) 
Moreover 
it seems that in general strict smoothness of a $\gt$ module $V$ and
$V(1) \in \GF$ do not imply that $V \in \Aff(\GF)_\ka$ (the problem
being that $V$ might have infinite length) but we do not know
a counterexample at this time.

As another consequence of the special case when $\f$ is a Cartan 
subalgebra of $\g$ (see part 4 of \exref{mainex})
one easily obtains that the fusion tensor product of 
a module in $\O_\ka$ and a module in the affine category $\O$ with 
central charge $\ka -h\spcheck$ is again a module in the affine category 
$\O.$ This was previously proved by Finkelberg \cite{F}. 
 
\thref{fin} opens up the possibility for defining translation functors 
\cite{Z, LW, Jan} in the categories $\Aff(\GF)_\ka$ using the 
Kazhdan--Lusztig 
fusion tensor product. It is also interesting to investigate if
tensoring with $\O_\ka$ can be used to define functors from the
categories $\Aff(\GF)_\ka$ to some categories of representations
of the corresponding quantum group $U_q(\g)$ in the spirit 
of Kazhdan and Lusztig \cite{KL4}. This can be viewed as a 
procedure of ``quantizing categories of modules over a complex 
simple Lie algebra'' by considering first categories of modules
for the related affine Kac--Moody algebra.

Acknowledgments: I am grateful to Dan Barbasch, Ivan Cherednik,
Michael Finkelberg, 
Edward Frenkel, David Kazhdan, and Gregg Zuckerman for very helpful 
discussions and correspondence. I also wish to thank Wolfgang Soergel 
for communicating a reference to us. 
\sectionnew{Properties of induced modules}
Recall that on any smooth $\gt$ module $V$ of central charge 
$\ka- h\spcheck$ $(\ka \neq 0)$ 
there is a well defined action of the Sugawara operators
(see \cite{Kac} for details)
\begin{equation}
\lb{Sug}
L_k = \frac{1}{2 \ka} \sum_{p} 
\sum_{j \in \Zset} \colon (x_p \e^{-j})(x_p \e^{j+k}) \colon
\end{equation}
where the first sum is over an orthonormal basis $\{ x_p \}$
of $\g$ with respect to the bilinear form $(.,.).$
In \eqref{Sug} the standard normal ordering is used,
prescribing pulling to the right the term $x \e^n$ with larger 
$n.$

The Sugawara operators define 
a representation of the Virasoro algebra on
$V$ \cite{Kac} for which
\begin{equation}
\lb{commut}
[L_k, x \e^n] = -n (x \e^{n+k}).
\end{equation}

For a $\gt$ module $V$ consider the generalized eigenspaces of 
the operator $L_0$
\[
V^\xi = \{ v \in V \mid (L_0 -\xi)^n v =0 \; 
\mbox{for some integer} \; n \}, \; \xi \in \Cset.
\]
Since $L_0$ commutes with $\g \hra \gt$ (see \eqref{commut})
each $V^\xi$ is a $\g$ module.

\bde{Weyl} For a $\g$ module $M$ and $\ka \in \Cset$ 
define the Weyl module
\[
\Ind(M)_\ka = U(\gt) \otimes_{U(\g[\e] \oplus \Cset \c)} M
\]
of $\gt$ where $\g[\e]$ acts through the quotient map
$\g[\e] \to \g[\e] / \e \g[\e] \cong \g$ and $\c$ acts by
$(\ka - h\spcheck) \id.$
\ede

\bpr{WeylSug} Assume that $M$ is a $\g$ module on which the 
Casimir $\Om$ of $\g$ acts by $a . \id$ for some $a \in \Cset$
and that $\ka$ is a nonzero complex number. Then:

(1) $\Ind(M)_\ka = \bigoplus_{\xi \in \Cset} \Ind(M)_\ka^\xi$ and 
$\Ind(M)_\ka^\xi$ are actual (not generalized) eigenspaces of 
$L_0.$

(2) $\Ind(M)_\ka^\xi= 0$ unless $\xi \in a/2\ka + \Zset_{\geq0}$
and as $\g$ modules
\begin{equation}
\lb{Sad}
\Ind(M)_\ka^{(a/2\ka +n)} \cong M \otimes S(\ad)^n
\end{equation}
where $S(\ad)^n$ denotes the
degree $n$ component of the symmetric algebra of the graded 
vector space $\g \op \g \op \dots$ with $k$-th term sitting 
in degree $k$ $(k=1,2, \ldots)$
considered as a $\g$ module under the adjoint action.

(3) If $M$ is an irreducible $\g$ module and $V$ is a nontrivial
$\gt$ submodule of 
$\Ind(M)_\ka$ then 
\[
V \cap \Ind(M)_\ka^{a/2\ka} = 
V^{a/2\ka} =0.
\]
\epr

Part 2 follows from \eqref{commut} and the 
Poincare--Birkhoff--Witt lemma. Parts 1 and 3 
are straightforward.

\bco{indirr} If $M$ is an irreducible $\gt$ module then 
$\Ind(M)_\ka$
has a unique maximal $\gt$ submodule $M_\max.$ It satisfies 
$M_\max \cap \Ind(M)_\ka^{a/2\ka} =0.$
The corresponding irreducible quotient will be denoted by
\[
\Irr(M)_\ka= \Ind(M)_\ka / M_\max.
\]
\eco

We now use a theorem of Kostant \cite{Kos}. Recall that through
the Harish-Chandra isomorphism the center $Z(\g)$ of $U(\g)$ is identified
with $S(\h)^W$ for a given Cartan subalgebra $\h$ of $\g,$ see
e.g. \cite[Chapter 7.4]{Dix}. Thus the characters of $Z(\g)$ 
are parametrized by $\h^\st/W.$ The character corresponding to
the $W$ orbit of $\la \in \h^\st$ will be denoted by 
$\chi_\la \colon Z(\g) \ra \Cset.$ Recall also that for the Casimir 
element $\Om \in Z(\g)$
\[
\chi_\la(\Om) = |\la|^2 - |\rho|^2
\]
where $\rho \in \h^\st$ is the half-sum of the positive roots of $\g$ for 
the Borel subalgebra used to define the Harish-Chandra isomorphism.

\bth{Kostant} (Kostant) Let $M$ be a $\g$-module with 
infinitesimal character $\chi_\la,$ $\la \in \h^\st$ and
$U$ be a finite dimensional $\g$ module with weights 
$\mu_1,$ $\cdots,$ $\mu_n,$ counted with their multiplicities. 
Then 
\[
\prod_{i=1}^n \left(z - \chi_{\la+\mu_i}(z). \id \right)
\]
annihilates $U \otimes M$ for all $z \in Z(\g),$ in particular 
\[
\prod_{i=1}^n \left( \Om -(|\la+ \mu_i|-|\rho|^2). \id \right)
\]
annihilates $U \otimes M.$
\eth  

\ble{submod} Let $M$ be an irreducible $\g$ module with infinitesimal 
character $\chi_\la,$ $\la \in \h^\st.$ Then for any two $\gt$ submodules 
$V$ and $V'$ of $\Ind(M)_\ka$ such that 
\[
V \sub V', \quad V \neq V'
\]
there exists an element $\mu$ of the root lattice 
$Q$ of $\g$ such that

(a) the operator $L_0: V'/V \ra V'/V$ has the eigenvalue 
\[
\frac{1}{2 \ka} \big( |\la +\mu|^2 - |\rho|^2 \big)
\]
and

(b) 
\[
\frac{1}{2 \ka} \big( |\la +\mu|^2 - |\la|^2 \big) \in \Zset_{\geq 0}.
\]
\ele

\begin{proof} 
The subspaces $V$ and $V'$ of $\Ind(M)_\ka$ are invariant under 
$L_0$ and therefore $L_0$ induces a well defined endomorphism
of $V/V'.$ Choose the eigenvalue $\xi_0$ of $L_0$ on $V/V'$ with 
minimal real part. (It exists due to part 2 of \prref{WeylSug}.)
Then $\gt_+$ annihilates $(V/V')^{\xi_0}$ because of \eqref{commut}. This 
implies that the Casimir of $\g$ acts on $(V/V')^{\xi_0}$ by 
$2 \ka \xi_0.\id.$ 
On the other hand $(V/V')^{\xi_0}$ considered as a $\g$ module is a 
subquotient of the module $\Ind(M)_\ka^{(|\la|^2- |\rho|^2)/2\ka +n}$
for some $n \in \Zset_{\geq 0},$ see 
\eqref{Sad} in \prref{WeylSug}.

According to Kostant's theorem 
\[
\xi_0 = \frac{1}{2\ka} \big( |\la + \mu|^2 - |\rho|^2 \big)
\]
for some $\mu$ in the root lattice $Q$ of $\g.$ At the same time
\[
\xi_0 = \frac{1}{2 \ka} \big( |\la|^2 - |\rho|^2 \big)  + n
\]
for the nonnegative integer $n$ above.
This weight $\mu$ satisfies properties (a) and (b) above.
\end{proof}

\leref{submod} motivates the following definition.
Assume that $M$ is a $\g$ module of finite length
and that all tensor products of $M$ with powers of the adjoint 
representation of $\g$ have finite length as well.
Then for every $\gt$ submodule $V$ of $\Ind(M)_\ka$
set
\begin{equation}
\lb{length}
\de(V)= \sum_{
     \begin{array}{c}
        \xi \colon 
        \xi \in (|\la|^2 - |\rho|^2)/2 \ka + \Zset_{\geq 0} \; 
        \mbox{and} \\
        \exists \mu \in Q, 
        \xi = (|\la+\mu|^2 - |\rho|^2)/2 \ka \\
     \end{array}        
             }
l(V \cap \Ind(M)_\ka^\xi)
\end{equation}
where $l(.)$ denotes the length of a $\g$ module. The $\g$ modules 
$\Ind(M)_\ka^\xi$ have finite lengths because of part 2 of 
\prref{WeylSug}. Clearly 
$\de(V) \in \Zset_{>0} \cup \{ \infty \}.$

The following lemma contains the
major property of the function $\de(.).$ 

\ble{decrease} Let $M$ be a $\g$ module for which the tensor products 
\[
M \otimes (\ad)^{\otimes n}
\]
have finite length for $n \in \Zset_{\geq 0}.$ If $V \sub V'$ are 
two $\gt$ submodules of $\Ind(M)_\ka$ such that $V \neq V'$ then 
either $\de(V) = \de(V')=\infty$ or $\de(V) < \de(V').$
\ele

In the case of an irreducible $\g$ module $M,$ 
\leref{decrease} is a direct consequence of \leref{submod},
analogously to the proof of \cite[Proposition 2.14]{KL}.
The general case follows from the exactness of the functor
$M \mt \Ind(M)_\ka.$

Next, for some $\g$ modules $M,$ we establish bonds on
$\de(V)$ for all $\gt$ submodules $V$ of 
$\Ind(M)_\ka.$ First note that 
\begin{equation}
\lb{lbound}
C = \min_{\mu \in Q} 
\Re \big(|\mu|^2 + 2(\la,\mu)
\big) 
\end{equation}
exists and is finite because
for a fixed $\la \in \h^\st,$ 
$\Re(|\mu|^2+ 2(\la,\mu))= |\mu|^2 + 2 \Re (\la, \mu)$ 
is a positive definite 
quadratic function on the root lattice $Q$ of $\g.$
Moreover $C\leq 0$ because the above function of
$\mu$ vanishes at $\mu=0.$

\ble{positive}
(1) If $\ka \notin \Rg$ then  
\begin{equation}
\lb{finite_l}
\de(\Ind(M)_\ka) < \infty.
\end{equation}

(2) Define the set 
\begin{equation}
\lb{Xla}
X_\la = \left\{ 
\frac{|\mu|^2 + 2(\la, \mu)}{2 n} 
\mid
\mu \in Q, n \in \Zset_{> 0} 
\right\} \sub \Cset.
\end{equation}
If $\ka \notin X_\la$ and in particular if 
\begin{equation}
\lb{Yla}
\ka \notin Y_\la = \Qset + 
\Qset (\la, \al_1) + \dots \Qset (\la, \al_r) 
\sup X_\la 
\end{equation}
then
\[
\de(\Ind(M)_\ka) = l(M).
\]

(3) If $\Re \ka < \frac{C}{2}$ then 
\[
\de(\Ind(M)_\ka) = l(M).
\]
\ele

\begin{proof} Since the sum \eqref{length} is 
over those $\xi \in \Cset$ for which
\[
\xi = (|\la|^2 - |\rho|^2)/2 \ka + n =  
(|\la+\mu|^2 - |\rho|^2)/2 \ka
\]  
for some $\mu \in Q,$ $n \in \Zset_{\geq 0},$ 
for each $\xi$ in \eqref{length} there exists
a pair $(\mu, n) \in Q \times \Zset_{\geq 0}$ such that
\begin{equation}
\lb{pair}
|\mu|^2 + 2 (\la, \mu) = 2 \ka n.
\end{equation}

(1) We claim that if $\ka \notin \Rg$ then the sum 
in \eqref{length} is finite for any $\gt$ submodule 
of $\Ind(M)_\ka.$
This follows from the fact
that if $\ka \notin \Rg,$ then the set of pairs \eqref{pair} 
is finite because for each $\epsilon >0$
there exists $R>0$ such that 
\[
\frac{\Im(|\mu|^2 + 2 (\la, \mu))}{\Re(|\mu|^2 + 2 (\la, \mu))}
<\epsilon\; \mbox{for} \; |\mu|>R.
\]
The statement now follows from the fact that for each of those finitely 
many $\xi$'s the $\g$ module $\Ind(M)^\xi_\ka$ has finite length.

(2) In this case the set of pairs \eqref{pair} consists only of
the pair $(\mu, n)=(0,0),$ i.e.
\[
\de(\Ind(M)_\ka) = l\big(\Ind(M)^{|\la|^2/2 \ka}_\ka \big)
= l(M).
\]

(3) If $\Re \ka < \frac{C}{2}$ then for every $\mu \in Q$ and 
$n \in \Zset_{>0}$ 
\[
\Re \ka < \frac{C}{2} \leq \frac{C}{2n} \leq 
\Re \frac{|\mu|^2 + 2(\la, \mu)}{2 n}
\]
because $C\leq 0,$ as noted before the statement of 
\leref{positive}. Therefore 
$\Re \ka < \frac{C}{2}$ implies that $\ka \notin X_\la$ 
and part (3) follows from part (2).   
\end{proof}

\bth{induced} Let $M$ be a $\g$ module as in \leref{decrease}.  

(1) If $\ka \notin \Rset_{\geq 0}$ 
then $\Ind(M)_\ka$ has finite composition series 
with quotients of the type $\Irr(M')_\ka$ for some irreducible 
subquotients $M'$ of $M \otimes S(\ad)^n,$ see \eqref{Sad}.

(2) If $\Re \ka < C/2$ and $M$ is irreducible then $\Ind(M)_\ka$ is an 
irreducible $\gt$ module.

(3) If $\ka \notin Y_\la$ or more generally $\ka \notin X_\la$ 
and $M$ is irreducible then
$\Ind(M)_\ka$ is again an irreducible $\gt$ module (see \eqref{Xla}, 
\eqref{Yla} for the definitions of the sets 
$X_\la \sub Y_\la  \sub \Cset$).
\eth

\noindent
\begin{proof} The first statement in part 1 and parts 2-3
follow from \leref{decrease} and \leref{positive}.

To prove the second statement in part 1, assume that
$V \sub V' \sub \Ind(M)_\ka$ are two submodules 
such that $V'/V$ is a nontrivial irreducible
$\gt$ module. Choose the eigenvalue $\xi_0$ of $L_0$ acting on $V'/V$ with 
minimal
real part. Then $(V'/V)^{\xi_0}$ 
is annihilated by $\gt_+$ and 
is an irreducible $\g$ module, otherwise 
if $M_0$ is a submodule of $(V'/V)^{\xi_0}$ 
we obtain a morphism $\Ind(M_0)_\ka \ra V'/V$
whose image is a nontrivial $\gt$ submodule of $V'/V$
because of \prref{WeylSug}.

Next we obtain a homomorphism 
$\Ind \big( (V'/V)^{\xi_0} \big)_\ka \ra V'/V$
which needs to be surjective and consequently we obtain that 
\[
V'/V \cong \Irr \big( (V'/V)^{\xi_0} \big)_\ka.
\]
\end{proof}

\bre{LZuck} Part 2 of \thref{induced} generalizes a result of
Lian and Zuckerman \cite[Proposition 2.2]{LZ2} in the case when $M$ is a 
Harish-Chandra module which they obtained, using the Jacquet functor.
\ere
\sectionnew{The categories $\Aff(\C)_\ka$}
{\em{Throughout this section we will assume that $\C$ is a full 
subcategory of the 
category of $\g$ modules of finite length which is closed under tensoring 
with the adjoint representation of $\g$ and taking subquotients,
see \exref{mainex}. We will also assume that}} 
\[
\ka \notin \Rset_{\geq 0}.
\] 

\bpr{ind} Under the above assumptions for any $M \in \C$ 

(1) $\Ind(M)_\ka$ has finite composition series with quotients of the type
$\Irr(M')_\ka$ for some subquotients $M'$ of $M \otimes S(\ad)^n,$
see \eqref{Sad}.

(2) $\Ind(M)_\ka \in \Aff(\C)_\ka.$
\epr

\begin{proof} Part 1 follows from the exactness of the 
functor $\Ind(.)_\ka$ and part 1 of \thref{induced}.

Since the functor $V \mt V(N)$ (from the category of $\gt$ modules to the
category $\g$ modules) is left exact to prove part 2 it is sufficient to
prove that 
\begin{equation}
\lb{relaxed}
\Irr(M)_\ka \in \Aff(\C)_\ka
\end{equation}
whenever $M \in \C$ is an irreducible $\g$ module. Indeed the left 
exactness shows that $\Ind(M)_\ka(N)$ is finitely generated by induction
on the length of $\Ind(M)_\ka.$ But $\Ind(M)_\ka(N)$ 
is also a submodule of 
$M\otimes \left( \bigoplus_{n \in \Zset_{\geq 0} } S(\ad)^n \right)$
and thus of the truncated tensor product for $n \leq k$ for
some integer $k.$
This implies that $\Ind(W)_\ka \in \C.$

To show \eqref{relaxed} we note that $\Irr(M)_\ka(1) \cong M$ for an 
irreducible $\g$ module $M.$ (Indeed $\Irr(M)_\ka(1)$ is an 
irreducible $\g$ module because if $M'$ is a submodule of it 
then there would exist a homomorphism of $\gt$ modules 
$\Ind(M')_\ka \ra \Irr(M)_\ka$ whose image would be
a nontrivial 
submodule of $\Irr(M)_\ka.$ Now \coref{indirr} gives 
$\Irr(M)_\ka(1) \cong M$.) 
Thus $\Irr(M)_\ka(1) \in \C.$
Finally we use the following result of
Kazhdan and Lusztig \cite[Lemma 1.10(d)]{KL}:

For any $\gt$ module $V$ there is an exact sequence of $\g$ modules
\[
0 \ra V(1) \ra V(N) \ra \Hom_\Cset(\g, V(N-1)) 
\cong \ad \otimes V(N-1), \; \mbox{for} \; N \geq 2
\]  
where the map $i$ is given by 
\[
i(v) (x)= (\e x) . v \in V(N-1), \; v \in V(N), \, x \in \g,
\]

By induction on $N$ one shows that $\Irr(M)_\ka(N)$ is
a finitely generated $\g$ module. Hence
$\Irr(M)_\ka(N)$ is a subquotient of 
$M \otimes \left( \bigoplus_{n=0}^k S(\ad)^n \right)$
for some sufficiently large integer $k$ and thus 
belongs to $\C.$
\end{proof} 

\bde{nil} A $\g[\e]$ module $\N$ is called a nil-$\C$-type 
module if $U(\gt_+)^n$ annihilates it
for a sufficiently large integer $n$
and considered as a $\g$ module $\N \in \C.$
\ede

\ble{nil2} A module $\N$ over $\g[\e]$ is a nil-$\C$-type module
if and only if it admits a filtration by $\g[\e]$ submodules
\[
\N=\N_m \sup N_{m-1} \sup \dots \N_1 \sup \N_0=0
\]
such that $\gt_+ \N_i \sub \N_{i-1}$ and $\N_i/\N_{i-1}$
are irreducible $\g$ modules which belong to $\C.$
\ele

For a $\g[\e]$ module $M$ we define the induced $\gt$ module
\begin{equation}
\label{I}
\I(M)_\ka = U(\gt) \otimes_{U(\g[\e] \oplus \Cset \c)} M
\end{equation}
where $M$ is extended to a $\g[\e] \oplus \Cset \c$ module
by letting $\c$ act by $\ka - h\spcheck.$

\bde{genW} A generalized Weyl module over $\gt$ of type $\C$ 
and central charge $\ka - h\spcheck$
is an induced module 
\[
\I(M)_\ka
\]
for some nil-$\C$-type module $\N$ over $\g[\e].$ 
\ede 

\bth{affineC} Let $V$ be a $\gt$ module of central charge 
$\ka-h\spcheck,$ $\ka \notin \Rset_{\geq 0}.$ 

(1)The following (a)-(c) are equivalent

(a) $V \in \Aff(\C)_\ka,$

(b) There exists a positive integer $N$ such that $V(N) \in \C$ and 
$V(N)$ generates V as a $\gt$ module,

(c) $V$ is a quotient of a generalized Weyl module of type $\C.$

(2) The irreducible objects in $\Aff(\C)_\ka$ are the modules 
$\Irr(M)_\ka$ for irreducible $\g$ modules $M.$ 
In addition for two nonisomorphic irreducible $\g$ 
modules $M$ and $M'$ the $\gt$ modules 
$\Irr(M)_\ka$ and $\Irr(M')_\ka$ are not isomorphic.

(3) The category $\Aff(\C)_\ka$ is closed under taking 
subquotients. Any module in $\Aff(\C)_\ka$ has finite length
and thus has a 
filtration with quotients of the type $\Irr(M)_\ka$
for some irreducible $\g$ modules $M.$
\eth

\begin{proof} Part 1: Obviously (a) implies (b).

Condition (b) implies (c) because assuming (b), $V(N)$ is a 
naturally a nil-$\C$-type module over $\g[\e]$ and thus $V$ is a 
quotient of the corresponding generalized Weyl module.

Condition (c) implies (a) as follows. Because of \leref{nil2} and
the exactness of the induction functor any generalized Weyl module
for $\gt$ has a filtration with quotients of the type $\Ind(M)$ 
for some irreducible $\g$ modules $M \in \C.$ 
Now \prref{ind} implies that it also 
has a filtration with quotients of the type $\Irr(M)_\ka$ 
(again for some irreducible $\g$ modules $M$).
Thus any quotient $V$ of a generalized Weyl module
has a filtration of the same type. The left exactness of the functor
$V \ra V(N)$ implies by induction that $V(N)$ are finitely generated 
$\g$ modules. Using \eqref{Sad} as in the proof
of \prref{ind} we see that $V(N) \in \C.$
  
Part 2: If $V$ is an irreducible $\gt$ module which belongs
to $\Aff(\C)_\ka$ then $V(1)$ should be an irreducible $\g$ module.
Otherwise, since it has finite length, it would contain 
an irreducible $\g$ module $M$ and one would obtain a 
homomorphism $\Ind(M)_\ka \ra V$ which should factor
through an isomorphism $\Irr(M)_\ka \cong V.$ But this 
is a contradiction since $\Irr(M)_\ka(1) \cong M,$ see the proof
of \prref{ind}.   

If $\Ind(M)_\ka$ and $\Irr(M')_\ka$ are isomorphic 
$\gt$ modules
for two irreducible $\g$ modules $M$ and $M'$ then
$\Ind(M)_\ka(1) \cong M$ and
$\Ind(M')_\ka(1) \cong M'$ are isomorphic
$\g$ modules and thus $M \cong M'.$

Part 3 follows from the characterization (c)
of $\Aff(\C)_\ka$ by generalized Weyl modules. 
It can be easily proved directly.
\end{proof}

As a consequence of part 1, condition (b) of \thref{affineC} one 
obtains: 
\bco{fingen+} Any module $V \in \Aff(\C)_\ka$ is 
finitely generated over $U(\g[\e^{-1}]).$
\eco

In the case when the category $\C$ of $\g$ modules is closed under 
extension we get that the category of $\gt$ modules is closed under 
extensions too. This is the case for the categories $\Fin_\g$ and
more generally $\GF$ in \exref{mainex} when $\f$ is a semisimple Lie 
algebra.

\bth{closedex} Assuming that the category $\C$ is closed under extensions
and $\ka \notin \Rset_{\geq 0}$ the following hold:

(1) The category of $\gt$ modules $\Aff(\C)_\ka$ is closed under 
extension inside the category of $\gt$ modules of 
central charge $\ka - h\spcheck$ and

(2) A $\gt$ module of central charge $\ka -h\spcheck$ belongs to 
$\Aff(\C)_\ka$ if and only if 
it has a finite composition series with quotients of the
type $\Irr(M)_\ka$ for some irreducible $\g$ modules $M \in \C.$
\eth

The proof of \thref{closedex} mimics the one of \thref{affineC}.

Finally we return to \prref{Sugchar}

\noindent
{\em{Proof of \prref{Sugchar}.}} If $V \in \Aff(\C)_\ka$ then it is a
quotient of a generalized Weyl module, \eqref{Sugdecom} and 
$V^\xi \in \C$ holds because of \prref{WeylSug}.

In the other direction -- assume that $V$ is a finitely generated $\gt$ 
module 
for which \eqref{Sugdecom} holds and $V^\xi \in \C.$ Then
$V$ is generated as a $U(\gt)$ module
by $v_j \in V^{\zeta_j}$ for some 
$\zeta_j \in \Cset,$ $j=1,\ldots, k.$ 
Thus it is generated by
\begin{equation}
\label{ssum}
\bigoplus_{\xi \colon 
\xi=\zeta_j - n, \; j=1,\ldots, k, \, n \in \Zset_{\geq 0} }
V^\xi.
\end{equation}
(Note that the above sum is finite because of \eqref{Sugdecom}).
Eq. \eqref{ssum} defines 
a $\g[\e]$ submodule of $V$ because of \eqref{commut} 
and as a $\g \hra \g[\e]$ module it belongs to
$\C$ since the sum in \eqref{ssum} is finite. Therefore $V$ is a quotient 
of the corresponding generalized Weyl module of type $\C$ and belongs to
$\Aff(\C)_\ka.$
\hfill \qed

There exists also a characterization of the categories
$\Aff(\C)_\ka$ in the spirit of Lian and Zuckerman \cite{LZ, LZ2}:

\bpr{LZchar} The category $\Aff(\C)_\ka$ consists exactly of those
$\gt$ modules $V$ of central charge $\ka -h\spcheck$ which are $\Zset$ 
graded
\begin{equation}
\label{Zgrad}
V = \bigoplus_{n \in \Zset} V_n
\end{equation}
with respect to the grading \eqref{grading} and
\[
V_n \in \C, \quad V_n =0 \; \mbox{for} \; n \ll 0.
\] 
\epr
\noindent
{\em{Sketch of the proof of \prref{LZchar}.}} 
Let $V \in \Aff(\C)_\ka.$ Then \eqref{Sugdecom} holds for some
$\xi_1,$ $\ldots,$ $\xi_n \in \Cset$ and we can 
assume that $\xi_i -\xi_j \notin \Zset.$ 
To get the grading \eqref{Zgrad}
we can set e.g.
\[
V_{m} := V_{\xi_1+m} \oplus \ldots \oplus V_{\xi_n+m}.
\]

The opposite statement is proved similarly to \prref{Sugchar}.
\hfill \qed \\  
\sectionnew{Duality in the categories $\Aff(\GF)_\ka$}
Let $\f$ be a subalgebra of $\g$ which is reductive in $\g$ and 
$\ka \notin \Rset_{\geq 0}.$ 
In this section, completely analogously to \cite{KL}, we define 
a natural duality in the categories $\Aff(\GF)_\ka.$
We will only state the results.

For any $\f$ module $M$ we define
\begin{equation}
\label{d}
M^d := (M^\st)^{\f - fin}
\end{equation}
where $(.)^\st$ stays for the full dual and $(.)^{\f -fin}$ denotes 
the $U(\f)$-finite part, i.e. the set all $\eta$ such that 
$\dim U(\f) \eta < \infty.$

It is well known that 
\[
M \mapsto M^d
\]
is an involutive antiequivalence of $\GF.$

Recall \cite{KL} that $\gt$ has the following
involutive automorphism
\begin{equation}
\label{sharp}
(x \e^k)^\sharp = x (-\e)^{-k}, \; k \in \Zset; \quad
(\c)^\sharp = -\c.
\end{equation}
For a $\gt$ module $V$ by $V^\sharp$ we will denote the twisting of
$V$ by this automorphism.

Recall also \cite{KL} that for any $\gt$ module $V$ the strictly smooth
part 
\[
V(\infty) = \cup_{N \in \Zset_{>0}} V(N)
\]
of $V$ is a $\gt$ submodule.

For any $\gt$ module $V$ define
\begin{equation}
\label{DD}
D(V) :=  (V^d)^\sharp(\infty) = \cup_{N \geq 1} (V^d)^\sharp(N).
\end{equation}
Here the restricted dual $V^d$ is defined with respect to the action $\f$ 
on $V$ coming from the embedding $\f \hra \g \hra \gt.$
It is clear that $D(V)$ is a strictly smooth $\gt$ module. If
$V$ has central charge $\ka -h\spcheck$ then $D(V)$ has the same 
central charge. 

It is easy to see that if $V \in \Aff(\GF)_\ka$ then the generalized 
eigenspaces $(V^d)^\xi$ of the Sugawara operator $L_0$ for the 
$\gt$ module $V^d$ are given by
\begin{equation}
\label{dualxi}
(V^d)^\xi =\{ \eta \in V^d \mid\eta(V^\zeta)=0 \; 
\mbox{for} \;\zeta \neq \xi \}.
\end{equation}

\bpr{Dual} (1) Fix $V \in \Aff(\GF)_\ka$ with decomposition 
\eqref{Sugdecom} for some $\xi_1,$ $\ldots,$ $\xi_n \in \Cset.$ Then as a 
subspace of $V^d$ the dual module $D(V)$ is 
\[
D(V) = \bigoplus_{\xi \colon \xi-\xi_1 \in \Zset_{\geq 0}, \ldots, 
\xi - \xi_n \in \Zset_{\geq 0} }(V^d)^\xi.
\]

(2) The contravariant functor $D$ is an involutive 
antiequivalence of the category $\Aff(\GF)_\ka.$
 
(3) The functor $D$ transforms simple objects 
$\Irr(M)_\ka \in \Aff(\GF)_\ka$ by
\[
D(\Irr(M)_\ka) \cong \Irr(M^d)_\ka.
\]
\epr

Parts 1 and 3 are proved analogously to Section 2.23 and
Proposition 2.24 in \cite{KL}. Similarly to
\cite[Proposition 2.25]{KL} 
one shows that the functor $D$ is exact.
This implies that for any 
$V \in \Aff(\GF)_\ka,$ $D(V)$ has finite length and thus belongs
to $\Aff(\GF)_\ka,$ e.g. because of \prref{Sugchar}. Now part
2 of \prref{Dual} is straightforward.
\sectionnew{Finiteness properties of the Kazhdan--Lusztig tensor product}
In this section we prove \thref{fin}. 

First we recall the definition of the Kazhdan--Lusztig fusion tensor 
product \cite{KL}.
Consider the Riemann sphere $\cp$ with three fixed distinct 
points $p_i,$ $i=0,1,2$ on it. Choose local coordinates (charts) at
each of them, i.e. isomorphisms $\ga_i: \cp \ra \cp$ such that
$\ga_i(p_i)= 0$ where the second copy of $\cp$ is equipped with 
a fixed coordinate function $\e$ vanishing at $0.$

Set $R=\Cset[\cp \backslash \{p_0, p_1, p_2\}]$
and denote by $\Ga$ the central extension
of the Lie algebra $\g \otimes R$
by
\begin{equation}
\label{Gaext}
[f_1 x_1,  f_2 x_2] := f_1 f_2 [x_1, x_2] +
\Res_{p_0} (f_2 d f_1) (x_1, x_2) \c,
\end{equation}
for $f_i \in R$ and $x_i \in \g.$ Here $(.,.)$ denotes the invariant
bilinear form on $\g,$ fixed in Section 1. There is a canonical
homomorphism 
\begin{equation}
\label{Gahom}
\Ga \ra \widehat{\g \oplus \g}, \quad
x f \mapsto (x \ev (\ga_1^\st)^{-1}(f), x \ev (\ga_2^\st)^{-1}(f)),
\; \c \mapsto -\c
\end{equation}
where $\ev(.)$ denotes the power series expansion of a rational 
function on $\cp$ at 0 in terms of the coordinate function $\e.$

Define 
\[
G_N = \Span \{ (f_1 x_1) \ldots (f_N x_N) \mid f_i 
\;\mbox{vanish at} \; p_0, \; x_i \in \g \} \sub U(\Ga).
\]

Fix two smooth $\gt$ modules $V_1$ and $V_2$ of central charge
$\ka - h\spcheck.$ 
Equip $W = V_1 \otimes_\Cset V_2$ with a structure of $\Ga$ module
with central charge $-\ka + h\spcheck$ using the homomorphism 
\eqref{Gahom}.
Clearly 
\[
W \supset G_1 W \supset G_2 W \supset \ldots
\]
and one can consider the projective limit of vector spaces
\begin{equation}
\label{proj}
\wh{W} = \lim_{\lla} W/G_N W.
\end{equation}
Define an action of $\gt$ on $\wh{W}$ as follows.
Fix $m \in \Zset$ and for each $n \in \Zset_{> 0}$ choose
$g_{n,m} \in R$ such that
\[
(\ga^\st_0)^{-1}(g_{n, m}) - \e^m \;
{\mbox{vanishes of order at least}} \; n \;
{\mbox{at}} \; 0.
\]
Set $\ol{m} = \max \{ -m, 0 \}$ and  
\begin{equation}
\label{gtact}
x \e^m . 
(w_1, w_2, \ldots) =
((x g_{1,m}) w_{\ol{m}+1}, (x g_{2,m}) w_{\ol{m}+1}, \ldots)
\end{equation}
for any sequence $(w_1, w_2, \ldots)$ in $W$ representing an element
of the projective limit \eqref{proj} i.e. $w_N \in W$ and 
$w_{N+1}-w_N \in G_N W.$ In \cite{KL} it is shown that this 
defines on $\wh{W}$ a structure $\gt$ module of central charge
$-\ka + h\spcheck,$ independent of the choice of $g_{n,m} \in \Ga.$ 
Finally the Kazhdan--Lusztig tensor product \cite{KL} of 
$V_1$ and $V_2$ is defined by 
\[
V_1 \dot{\otimes} V_2 := (\wh{W})^\sharp(\infty).
\]

We show finiteness properties of a dual construction of the fusion 
tensor product.  
Let $f_0$ be a rational function on $\cp$ (unique up to a multiplication
by a nonzero complex number)
having only one (simple) zero
at $p_0$ and only one (simple) pole at $p_1.$ 
For instance when $\ga_0(p_1)$ is finite 
$f_0(\e) = a \ga^\st_0(\e)/(\ga^\st_0(\e)-\ga_0(p_1)),$ $a \neq 0.$ Set
\begin{equation}
\label{XN}
X_N = \Span \{ (f_0 x_1) \ldots (f_0 x_N) \mid x_i \in \g \} \sub U(\Ga).
\end{equation}
Clearly $X_N \sub G_N.$

Kazhdan and Lusztig proved the following Lemma.

\ble{tensorle}\cite[Proposition 7.4]{KL} Assume that $V_i$ are two
strictly smooth $\gt$ modules of central charge $\ka-h\spcheck,$
generated by $V_i(N_i),$ respectively. Then
\[
V_1 \otimes V_2 = \sum_{k=0}^{N-1} X_k (V_1(N_1) \otimes V_2(N_2))
+ G_N (V_1 \otimes V_2)
\]
for all $N \in \Zset_{>0}.$
\ele

For a given category of $\g$ modules $\C$ denote by $\ol{\Aff}(\C)_\ka$ 
the category of strictly smooth $\gt$ modules $V$ such that
\[
V(N) \in \C.
\]
It differs from the category $\Aff(\C)_\ka$ in that we drop the 
condition for finite length.

Let $\f$ be a subalgebra of $\g$ which is reductive in $\g.$ Consider  
a module $U \in \O_\ka$ and a module $V \in \Aff(\GF)_\ka.$
Using the homomorphism \eqref{Gahom} $W= U \otimes_\Cset V$ becomes 
a $\Ga$ module of central charge $-\ka + h\spcheck.$
Note that the restricted dual $(U \otimes V)^d$ 
(recall \eqref{d}) is naturally a $\Ga$ submodule
of the full dual to $U \otimes V,$ both of central charge
$\ka - h\spcheck.$ (The restricted dual is taken with respect to
the embedding $\f \hra \g \hra \Ga$ using constant functions
on $\cp.)$
Following \cite{KL} define the following $\Ga$ submodule
of $(U \otimes V)^d$ 
\begin{equation}
\label{T}
T'(U, V) := \cup_{N \geq 1} T(U, V)\{N\} 
\end{equation}
where
\[
T(U, V)\{N\} := \Ann_{G_N} (U \otimes V)^d.
\] 
Eq. \eqref{T} indeed defines a $\Ga$ submodule of $(U \otimes V)^d$ 
since for any
$y \in \Ga$ and any integer $N$ there exists an integer $i$ such that
$G_{N+i} y \in \Ga G_N.$ In other words $T'(U,V)\{N\}$ is defined as 
\[
T'(U,V)\{N\} = \{ \eta \in (U \otimes V)^\st \mid \eta(G_N W) = 0, \; 
\dim U(\f) \eta < \infty \}.
\]
Similarly to \cite[Section 6.3]{KL} $T'(U,V)$ has a canonical action of 
$\gt$ 
(``the copy attached to $p_0$'') defined as follows. 
Let $\eta \in T(U,V)\{N\}.$ Fix $\om \in \Cset[\e, \e^{-1}],$ $x \in \g$ 
and 
choose $f \in R$ such that $f-\ga_0^\st(\om)$ has a zero of order at least 
$N$ at $p_0.$ Then
\begin{equation}
\label{tens}
(\om x) \eta := (f x) \eta
\end{equation}
correctly defines a structure of smooth $\gt$ module on $T'(U,V).$

\ble{functor} (1) For any two modules $U \in \O_\ka$ and $V \in 
\Aff(\GF)_\ka$, $T'(U,V)$ is a strictly smooth $\gt$ module 
of central charge $\ka - h\spcheck$ and
\begin{equation}
\label{nth}
T'(U,V)(N) = T'(U,V)\{N\} \in \GF
\end{equation}
considered as a $\g$ module. Thus 
\[
T' \colon \O_\ka \times \Aff(\GF)_\ka \ra \ol{\Aff}(\GF)_\ka, 
\quad
U \times V \mapsto T'(U,V)
\]
is a contravariant bifunctor.

(2) The bifunctor $T'$ is right exact in each argument.
\ele

The proof of \eqref{nth} is straightforward as  
\cite[Lemma 6.5]{KL} for the case $\f=\g$ $(\GF=\O_\ka).$
The part that $T'(U,V) \in \GF$ as a $\g$ module follows from
\leref{tensorle} and the Kostant theorem \cite{Kos} 
that $\GF$ is closed under tensoring 
with finite dimensional $\g$ modules (and taking subquotients).

Part 2, as in \cite[Proof of Proposition 28.1]{KL4}, follows from the left 
exactness of the functors $M \mapsto Ann_A(M)$ and 
$M \mapsto M^{A-fin}$ for a given algebra $A$ on the category of
all $A$ modules. (Here $M^{A-fin}$ denotes the $A$ finite part of $M,$
i.e. the space of $m \in M$ such that $\dim A. m < \infty,$
cf. section 4.)

\thref{fin} would be derived from the following Proposition.

\bth{fin2} For any modules $U \in \O_\ka$ and $V \in \Aff(\GF)_\ka$
the $\gt$-module $T'(U,V)$ has finite length, i.e.
\[
T'(U,V) \in \Aff(\GF)_\ka.
\]
\eth

Due to the right exactness of the bifunctor $T'$ it is sufficient to 
prove the following Lemma. {\em{This is the main part of the proof of 
\thref{fin}.}}

\ble{fin3} For any irreducible $\g$ modules $U_0 \in \Fin_\g$ and
$V_0 \in \GF$ the $\gt$ module $T'(\Irr(U_0)_\ka, \Irr(V_0)_\ka)$
has finite length.
\ele

To prove \leref{fin3} we can further restrict ourselves to charts 
$\ga_0 \colon \cp \ra \cp$ around $p_0$ such that
\begin{equation}
\label{ga0}
\ga_0(p_1) = \infty.
\end{equation}
This follows from the simple Lemma:
\ble{aux} 
Assume that the $\gt$ module $T'(U,V)$ associated to one chart
$\ga_0 \colon \cp \ra \cp$ around $p_0$ has finite length. Then
the respective module associated to any other chart 
$\al_0 \colon \cp \ra \cp$ around $p_0$ has finite length too.
\ele
\begin{proof}
Notice that $G_N$ does not depend on the choice of chart $\ga_0$ or 
$\al_0$ around $p_0$ and consequently the $\Ga$ module $T'(U,V)$ does not 
depend on such a choice either. Denote 
the actions of $\gt$ on the space $T'(U,V)$ associated to the charts
$\ga_0$ and $\al_0$ by $\mu$ and $\nu,$ respectively. We claim
that
\begin{equation}
\label{equa}
\mu(U(\gt)) \eta = \nu(U(\gt)) \eta \quad
\mbox{for any} \quad 
\eta \in T'(U,V).
\end{equation}
Since $\al_0 (\ga_0)^{-1}$ is an automorphism of $\cp$ that preserves the 
origin
\[
(\ga_0^\st)^{-1} \al_0^\st (\e) = \frac{a \e}{b \e +d}
\]
for some complex numbers $a \neq 0,$ $b,$ and $d \neq 0.$
Let $\eta \in T'(U,V)\{N\}.$ Then for each $n \in \Zset,$ $n \leq N$ 
there exist complex numbers $b_n,$ $\ldots,$ $b_N$ such that
\[
\al_0^\st (\e^n) - \sum_{k =n}^N b_k \ga_0^\st(\e^k) \, 
\mbox{vanishes of order at least} \, N \, \mbox{at} \, p_0.
\] 
This implies that
\[
\nu(x \e^n)\eta = \sum_{k =n}^N b_k \mu(x \e^k) \quad 
\mbox{for any} \quad x \in \g
\]
and by induction \eqref{equa}.

Assuming that the $\gt$ module $(T'(U,V), \mu)$ has finite length we get
that $(T'(U,V), \nu)$ is finitely generated $U(\gt)$ module which belongs
to $\ol{\Aff}(\GF)_\ka$ and thus it also has finite length due to
part 3 of \thref{affineC}.
\end{proof}

In the rest of this section we prove \leref{fin3} for a chart
$\ga_0$ around $p_0$ with the property \eqref{ga0}. 
Let us fix such a chart. Then 
\[
\ga_0^\st(\e^k) \in R = \Cset[\cp \backslash \{p_0, p_1, p_2 \} ], \quad
\mbox{for all} \quad k \in \Zset.
\]
As a consequence of this there exists an embedding
\begin{equation}
\label{morph}
\gt \hra \Ga, \quad x \e^k \mapsto x \ga^\st_0 (\e^k), \; \c \ra \c.
\end{equation}
In addition in \eqref{XN} the function $f_0$ can be taken simply as
$\ga_0^\st(\e).$

Fix two modules $U_0 \in \Fin_\g$ and $V_0 \in \GF$ and denote by
\[
W = \Irr(U_0)_\ka \otimes \Irr(V_0)_\ka 
\]
the related $\Ga$ module. Using the homomorphism \eqref{morph}
it becomes a $\gt$ module of central charge $-\ka + h\spcheck.$
We will denote by $W^\sharp$ the twisting of this $\gt$ module by 
the automorphism
$(.)^\sharp$ of $\gt$ (see \eqref{sharp}).
Note that $W^\sharp$ has central charge $\ka - h\spcheck.$ 

We claim that $\g[\e] \hra \gt$ preserves 
\[
W_0^\sharp := \Irr(U_0)_\ka(0) \otimes \Irr(V_0)_\ka(0) \sub W^\sharp.
\]
This follows from the facts that $x \e^n \in \gt$ acts on $W_0^\sharp$ by 
\[
x \ev (\ga_1^\st)^{-1} \ga_0^\st (-\e)^{-n} \otimes \id +
\id \otimes x \ev (\ga_1^\st)^{-1} \ga_0^\st (-\e)^{-n}
\]
(see \eqref{Gahom}), that $\ga_0^\st(\e^{-n})$ are regular functions on 
$\cp \backslash \{ p_0 \}$ for $n \in \Zset_{\geq0},$ 
and that $\g[\e] \hra \gt$ preserves 
$\Irr(U_0)_\ka(0)$ and $\Irr(V_0)_\ka(0).$

Consider the canonical induced homomorphism of $\gt$ modules
\[
\rho \colon \I(W_0^\sharp)_\ka \ra W^\sharp
\]
(Recall the definition \eqref{I} of an induced $\gt$ module $\I(.)_\ka$.) 
Dually we obtain a homomorphism of $\gt$ modules
\begin{equation}
\label{dualrho}
\rho^{d \sharp}: W^d \cong ((W^\sharp)^d)^\sharp \ra 
((\I(W_0^\sharp)_\ka)^d)^\sharp,
\end{equation}
see \eqref{d}.

Notice that in the definition \eqref{tens} of the action of $\gt$ on
$T'(U,V)$ the function $f$ can be simply taken as $\ga_0^\st(\om)$
because all functions $\ga_0^\st(\e^k),$ $k \in \Zset$ are regular
outside $\{ p_0, p_1 \}.$ This means that the structure of $\gt$
module on the space $T'(\Ind(U_0)_\ka, \Ind(V_0)_\ka)$ is simply the one 
induced from the $\Ga$ action by the homomorphism \eqref{morph}.
Thus $T'(\Ind(U_0)_\ka, \Ind(V_0)_\ka)$ is naturally a $\gt$ submodule
of $W^d.$

\ble{incl} The homomorphism $\rho^{d \sharp}$ \eqref{dualrho} restricts to 
an 
inclusion
\[
\rho^{d \sharp} \colon T'(\Ind(U_0)_\ka, \Ind(V_0)_\ka) \hra 
D(\I(W_0^\sharp)_\ka).
\]
(Recall that $D(\I(W_0^\sharp)_\ka)$ is the smooth part of 
$((\I(W_0^\sharp)_\ka)^d)^\sharp.$)
\ele
\begin{proof} Since 
$T'(\Ind(U_0)_\ka, \Ind(V_0)_\ka)\{N\} \sub W^d(N)$
\[
\rho^{d \sharp}(T'(\Ind(U_0)_\ka, \Ind(V_0)_\ka) \{N\}) \sub
((\I(W_0^\sharp)_\ka)^d)^\sharp(N)
\]
and
\[
\rho^{d \sharp}(T'(\Ind(U_0)_\ka, \Ind(V_0)_\ka) ) \sub
((\I(W_0^\sharp)_\ka)^d)^\sharp(\infty) =
D(\I(W_0^\sharp)_\ka).
\]

To show that this restricted $\rho^{d \sharp}$ is an inclusion assume that 
\\ $\eta \in T'(\Ind(U_0)_\ka, \Ind(V_0)_\ka)\{N\}$ is such that
$\rho^{d \sharp}(\eta)=0.$ Then
\begin{align}
\nn
(-1)^N \eta( \pi_1(x_1 \ga_0^\st(\e)) \ldots \pi_1(x_n \ga_0^\st(\e)) w_0) 
&= \eta( \pi_2(x_1 \e^{-1}) \ldots \pi_2(x_n \e^{-1}) w_0)
\\
\nn
&=
\rho^{d \sharp}
(\eta)( \pi_3(x_1 \e^{-1}) \ldots \pi_3(x_n \e^{-1}) w_0) = 0
\end{align}
for all $x_i \in \g,$ $w_0 \in W_0,$ $n \in \Zset_{\geq0}.$ Here
$\pi_1$ denotes the actions of $\Ga$ on $W,$ and
$\pi_2,$ $\pi_3$ denote the actions of
$\gt$ on $W^\sharp,$ $I(W^\sharp_0)_\ka,$ respectively. 
This means that 
\[
\eta|_{X_n W_0} =0, \; n \geq 0 \quad \mbox{and} \quad 
\eta|_{G_N W} =0. 
\]
Because of \leref{tensorle} $\eta=0.$
\end{proof}

Recall the canonical isomorphisms
\[
\Irr(U_0)_\ka(0) \cong U_0, \; 
\Irr(V_0)_\ka(0) \cong V_0.
\]
Kostant's theorem \cite[Theorem 3.5]{Kos} implies that $U_0 \otimes V_0$ 
has finite length as a $\g$ module, and thus
$W_0^\sharp = \Irr(U_0)_\ka(0) \otimes \Irr(V_0)_\ka(0)$
is a finite length $\g$ module.

Now \leref{fin3} follows from the following fact.

\bpr{finind} Let $M$ be a $\g[\e]$ module which is of finite length 
over $\g \hra \g[\e].$ Then the $\gt$ module 
$D(\I(M)_\ka)$ has finite length for 
$\ka \notin \Rset_{\geq 0}.$
\epr

\begin{proof} Denote by $\pi$ the action of $\g[\e]$ on $M.$ 
Define a new action $\ol{\pi}$ of $\g[\e]$ on the same space by
\[
\ol{\pi}(x \e^n) = \delta_{n, 0} \pi(x).
\] 
(This defines a representation of $\g[\e]$ since $\e \g [ \e ]$ is 
an ideal of $\g [ \e ].)$
This representation will be denoted by $\ol{M}.$ The underlining vector 
spaces of $M$ and $\ol{M}$ will be always identified. 

Consider the two $\gt$ modules $\I(M)_\ka$ and $\I(\ol{M})_\ka$
and identify their underlining spaces with 
\[
\M = U(\e^{-1}\g[\e^{-1}]) \otimes_\Cset M.
\] 
They are isomorphic as $\g[\e^{-1}]$  
modules and are naturally graded as $\g[\e^{-1}]$ modules with respect to 
the grading \eqref{grading} by
\[
\deg u \otimes m = -k \quad \mbox{for} \quad 
u \in U( \e^{-1} \g [ \e^{-1} ])^{-k}, m \in M.
\]
By $(.)^k$ we denote the $k$-th graded component of
a graded vector space (algebra). 

Then 
\[
U(\g[\e])^{-k} \I(M)_\ka =
U(\g[\e])^{-k} \I(\ol{M})_\ka=
\sum_{j=k}^\infty \M^{-j}.
\]
Denote 
\[
(\M^d)^{-k} = \{ \eta \in \M^d \mid \eta(\M^{-j}) = 0 \; 
\mbox{for} \; j \neq k \}.
\] 
In the definition of the restricted dual above, recall \eqref{d}, 
we use the $\f$ module structure on $\M$ coming from the 
identification of $\M$ with the isomorphic $\g$ modules
$\I(M)_\ka$ and $\I(M)_\ka.$ In other
words as an $\f \hra \g$ module $\M$ is the tensor product of
$U(\e^{-1}\g[\e^{-1}])$ (under the adjoint action) and $M$
(equipped with either the action $\pi$ or $\ol{\pi}$ which coincide 
when restricted to $\g).$

As subspaces of $\M^d$ 
\begin{equation}
\label{DN}
D( I(M)_\ka)(N) = D( I(\ol{M})_\ka)(N) = 
\sum_{j=0}^{N-1} (\M^d)^{-j}.
\end{equation}
This implies that the representation spaces of 
$D( I(M)_\ka)$ and $D(I(\ol{M})_\ka)$ 
can be identified with
\[
\bigoplus_{j=0}^\infty (\M^d)^{-j}.
\]
The actions of $\gt$ on this vector space related to 
$D(I(M)_\ka)$ and $D(I(\ol{M})_\ka)$ will be denoted 
by $\sig^\st$ and $\ol{\sig}^\st.$ 

We claim that: 
\begin{equation}
\label{good}
\mbox{If} \; \eta \in (\M^d)^{-j} \; \mbox{and} \;
g \in U(\g[\e^{-1}])^{-k} 
\; \mbox{then} \;
\sig^\st(g) \eta - \ol{\sig}^\st(g) \eta \in 
\bigoplus_{i=0}^{j+k-1} (\M^d)^{-i}.
\end{equation}
It suffices to check \eqref{good} for $g = x \e^{-k}.$ For this we need 
to show that if $u \in U(\e^{-1}\g[\e^{-1}])^{-p}$ and $p \geq k+j$
then
\begin{equation}
\label{check}
\big( \sig^\st(x \e^{-k}) \eta - 
\ol{\sig}^\st(x \e^{-k}) \eta \big)
(u \otimes m) =0.
\end{equation}
Let
\begin{equation}
\label{vp}
(-1)^k (x\e^k) u = \sum_i a_i b_i c_i
\end{equation}
for some $a_i \in U(\e^{-1} \g[\e^{-1}])^{-p+k-\vp(i)},$
$b_i \in U(\g),$ $c_i=1$ if $\vp(i)=0$ and 
$c_i \in U(\e \g[\e])^{\vp(i)}$ if $\vp(i)>0.$
Here $\vp$ is a map from the index set in the RHS of \eqref{vp} to
$\Zset_{\geq 0}.$ Then
\begin{align}
\nn
\big( \sig^\st(x \e^{-1}) \eta \big) (u \otimes m) &=
\sum_i \eta (a_i \otimes \sig(b_i c_i) m) \\
\nn
&= \delta_{p, k+j} 
\sum_{i\, : \vp(i)=0}  \eta (a_i \otimes \pi(b_i) m).
\end{align}
The second equality follows from $\eta \in (\M^d)^{-j}$ and
$p - k + \vp(i) >j$ unless $p=k+j$ and $\vp(i)=0.$

The same formula holds for $\ol{\sig}^\st$ with $\pi$ substituted in the 
RHS by $\ol{\pi}.$
The compatibility of $\pi$ and $\ol{\pi}$ on 
$\g \hra \g[\e]$ implies \eqref{check}.

According to 
\prref{Dual} $D(\I(\ol{M})_\ka) \in \Aff(\GF)_\ka.$ From 
\coref{fingen+} we get that $D(\I(\ol{M})_\ka)$ is finitely
generated as a $U(\g[\e^{-1}])$ module. We can assume that
it is generated by some homogeneous elements
\[
\eta_i \in (\M^d)^{-j_i}, \; i=1, \ldots, n.
\]
Then by induction \eqref{good} easily gives that 
$D(\I(M)_\ka)$ is generated
as a $U(\g[\e^{-1}])$ module by the same set 
$\{ \eta_1, \ldots, \eta_n \}.$ Thus $D(\I(M)_\ka)$ is 
finitely generated as a $U(\g)$ module and 
\[
D(\I(M)_\ka)(N) \in \GF
\]
because of \eqref{DN} which shows that
$D(\I(M)_\ka) \in \Aff(\GF)_\ka.$
\end{proof} 

Now as in \cite{KL} \thref{fin2} easily implies \thref{fin} and the 
following
\bpr{iso} In the setting of \thref{fin} and \thref{fin2} if
$U \in \O_\ka$ and $V \in \Aff(\GF)_\ka$ then
$U \dot{\otimes} V$ and $D(T'(U,V))$ are naturally isomorphic.
\epr

\end{document}